\documentclass{amsart}

\usepackage{psfig,latexsym,amssymb}
\def\proof{\medskip\noindent{\sc Proof. }}
\def\EOP{\hfill$\Box$}

\def\complex{\Bbb C}
\def\bdry{\partial }
\def\plus{\bigoplus }
\def\tensor{\bigotimes }

\newtheorem{thm}{Theorem}[section]

\newtheorem{defn}{Definition}[section]

\begin{document}
\title[The trivial representation in 
rank-selected homology]{Multiplicity of the trivial representation in 
rank-selected homology of the partition lattice}
\author{Phil Hanlon}
\curraddr{
        Phil Hanlon,
        Department of Mathematics,
        University of Michigan,
        Ann Arbor, MI 48109}
\email{hanlon@umich.edu}
\author{Patricia Hersh}
\curraddr{
        Patricia Hersh,
        Department of Mathematics, 
        University of Michigan,
        Ann Arbor, MI 48109}        
\email{plhersh@umich.edu}
\thanks{This work was supported by an NSF-AWM Mentoring Travel Grant.
The first author was also supported by NSF Grant DMS-9500979.}
\maketitle

\begin{abstract}
We study the multiplicity $b_S(n)$ of the trivial representation in the
symmetric group representations $\beta_S$ on the (top) homology of the 
rank-selected partition lattice $\Pi_n^S$.  We break the possible rank
sets $S$ into three cases: (1) $1\not\in S$, (2) $S=1,\dots ,i$ for 
$i\ge 1$ and 
(3) $S=1,\dots ,i,j_1,\dots ,j_l$ for $i,l\ge 1$, $j_1 > i+1$.  
It was previously shown by Hanlon 
that $b_S(n)=0$ for $S=1,\dots ,i$.  We use
a partitioning for $\Delta(\Pi_n)/S_n$ due to Hersh to 
confirm a conjecture of Sundaram [Su] that $b_S(n)>0$ for $1\not\in S$.
On the other hand, we use the spectral sequence of a filtered complex to show
$b_S(n)=0$ for $S=1,\dots ,i,j_1,\dots ,j_l$ unless a certain type of chain
of support $S$ exists.  The partitioning for $\Delta(\Pi_n)/S_n$ allows us
then to show that a large class of
rank sets $S=1,\dots ,i,j_1,\dots ,j_l$ for which 
such a chain exists do satisfy $b_S(n)>0$.
We also generalize the partitioning for $\Delta(\Pi_n)/S_n$ to 
$\Delta(\Pi_n)/S_{\lambda }$; when
$\lambda = (n-1,1)$, this partitioning leads to a proof of a conjecture
of Sundaram about $S_1\times S_{n-1}$-representations on the homology of 
the partition lattice.
\end{abstract}

\section{Introduction}

The natural action of the symmetric group $S_n$ on $\{ 1,\dots ,n\} $ 
gives rise to a rank-preserving, order-preserving action on the lattice 
$\Pi_n$ of partitions of $\{ 1,\dots ,
n\} $ ordered by refinement.  The resulting $S_n$-action permuting chains of 
comparable poset elements 
yields an $S_n$-representation on the top homology of the order
complex of the partition lattice.  We study the 
multiplicity $b_S(n)$ of the trivial representation in the representation
$\beta_S $ of the symmetric group $S_n$ on the homology of the partition 
lattice $\Pi_n $ restricted to
rank set $S$ for various $S\subseteq [n-2]$.  Questions about
these multiplicities were first suggested in [St1] and studied quite 
extensively, using symmetric functions, in [Su].  

We approach these 
questions from two other angles: spectral sequences of filtered complexes
and partitioning of quotient complexes.  A
partitioning of the quotient complex $\Delta (\Pi_n)/S_n$ 
lends itself well to proving lower bounds on $b_S(n)$, 
while spectral sequences of filtered complexes
seem well-suited to giving upper bounds.  One
of our interests is finding cases where we can get the two bounds to meet
and seeing how the two very different methods make use of the same 
information.  In particular, we give results
about when $b_S(n)$ is positive 
and when it is 0 (as well as when a related multiplicity $b_S'(n)$ is 
positive), including proofs of two conjectures of Sundaram [Su].

Recall that the {\bf order complex} 
$\Delta (P)$ of a finite poset $P$ with minimal and maximal elements
$\hat{0} $ and $\hat{1}$ is the simplicial complex comprised of an 
$i$-face for each chain $\hat{0} < u_0 < \cdots < u_i < \hat{1} $ of
comparable poset elements. 
Whenever a group $G$ acts on a graded poset $P$ in a rank-preserving,
order-preserving fashion, the group also permutes the poset chains, or
equivalently the faces in its order complex.  
This $G$-action on $\Delta (P)$ commutes with the 
boundary map, so the action on chains also gives rise to a 
$G$-representation on each of the homology groups of $\Delta (P)$.  
The action on chains also may be restricted to any rank
set $S$ giving rise to a group representation $\alpha_S$ permuting 
the chains of support $S$ and to representations on the homology of 
$\Delta (P)$ restricted to rank set $S$.

We will be interested in the alternating
sum $\beta_S = \sum_{T\subseteq S} (-1)^{|S-T|}\alpha_T $ 
of $S_n$-representations $\alpha_T $ on chains.
When $P$ is a Cohen-Macaulay poset, then $\Delta (P)$ (and each of its
rank-selected subcomplexes) only has top 
homology, in which case $\beta_S $ is the $G$-representation on the top
homology group in the rank-selected complex $\Delta^S $ obtained by 
restricting $\Delta $ to rank set $S$.  
The partition lattice is a Cohen-Macaulay poset, and we will
be interested in the multiplicity $b_S(n)$ of the trivial representation
in $\beta_S(n)$.  Our results about when $b_S(n)$ is positive come out of an 
analysis of the flag
$h$-vectors of the quotient complex $\Delta(\Pi_n)/S_n$, defined as follows:

\begin{defn} The {\bf flag $f$-vector} of a $(d-1)$-dimensional balanced
simplicial complex $\Delta $ is a vector with coordinates $f_S$ for each
subset $S\subseteq \{ 1,\dots ,d\} $ 
of the set of vertex colors for $\Delta $ where $f_S$ counts how many faces
in $\Delta $ have vertices colored exactly by $S$.  The {\bf flag $h$-vector}
is the alternating sum $h_S = \sum_{T\subseteq S} (-1)^{|S|-|T|}f_T$, or
alternatively, $h_S = (-1)^{|S|}\tilde{\chi }(\Delta^S )$.  See [St2] for
more background.
\end{defn}

\medskip 
In our case, the vertices of $\Delta (\Pi_n)$ are colored by poset rank,
and $\Delta (\Pi_n)$ is balanced because no two elements of a chain have
the same rank, implying no two vertices in the same face in $\Delta (\Pi_n)$
are assigned the same color.

The quotient complex $\Delta(P)/G$ consists of the $G$-orbits of 
faces in $\Delta (P)$, and it inherits the balancing by poset rank
from $\Delta (P)$ when $P$ is graded 
and $G$ preserves rank.  Note that $\Delta(P)/G$
typically is not the same as the order complex $\Delta(P/G)$ of the quotient
poset, and in particular $\Delta (\Pi_n/S_n) \ne \Delta (\Pi_n)/S_n$; 
there are elements $u<v, u'<v' \in P$ such that $u=\sigma u',
v=\tau v'$ for $\sigma ,\tau\in G$ but $u<v$ is not in the same 
$G$-orbit as $u'<v'$.  The quotient complex often is not a simplicial
complex, but it is always a boolean cell complex, i.e. a regular cell 
complex in which each cell has the combinatorial type of a simplex.

The multiplicity
$\langle \alpha_S ,1\rangle $ of the trivial representation within the 
group action $\alpha_S$ on chains of support $S$ equals the number of 
orbits in the action $\alpha_S$, 
i.e. it equals $f_S (\Delta(P)/G)$.  As observed in [Re], this implies that
$$\langle \beta_S ,1\rangle = \bigg\langle \sum_{T\subseteq S} (-1)^{|S-T|}
\alpha_T , 1\bigg\rangle = \sum_{T\subseteq S} (-1)^{|S-T|} f_T(\Delta(P)/G)
= h_S (\Delta(P)/G).$$  Hence, we 
will study flag $h$-vectors of quotient complexes
as a way of getting at $b_S(n)=\langle \beta_S ,1\rangle $.  
In particular, we will use the fact that
when a balanced complex $\Delta $ is shellable or partitionable, then 
$h_S(\Delta )$ counts minimal faces of support $S$ in the shelling or 
partitioning.  

\begin{defn} A pure simplicial complex $\Delta $ is {\bf partitionable} if the 
set of faces may be partitioned into a direct sum 
$$\Delta = [G_1,F_1]\cup\cdots\cup [G_k,F_k]$$ of intervals of 
boolean type where 
$F_1,\dots ,F_k$ are the facets of $\Delta $ and $G_i$ is a face of $F_i$ 
for $1\le i\le k$.  The complex $\Delta $ is {\bf shellable} if the facets may 
be ordered $F_1,\dots ,F_k$
so that for $2\le j\le k$,
the set $F_j\setminus \cup_{i<j} F_i$ of faces belonging to $F_j$ 
but not to any earlier facet, has a unique minimal element $G_j$.  
Thus, a shelling may be viewed as a type of partitioning.
\end{defn}

\medskip
Further background may be found in [St2].
We will use a very complicated 
partitioning for $\Delta(\Pi_n)/S_n$ given in [He] to show $b_S(n)>0$ for
various classes of $S$ by exhibiting facets $F_i$ with minimal faces $G_i$ 
of support $S$.  The partitioning for $\Delta(\Pi_n)/S_n$ has the property
that for a very large class of facets $F_i$, the minimal faces $G_i$ may
be described in terms of a generalized notion of ascents and descents in a
chain-labeling on orbits of saturated chains in $\Pi_n^* $.
Our strategy is to construct
facets achieving various descent sets $S$ to show that $b_S(n)>0$ for these
rank sets $S$.  

Denote by $\Pi_n^S$ the rank-selected subposet of the partition lattice
consisting of those poset elements of rank $r$ for some $r\in S\subseteq
[n-2]$.
We show in Section 3 that $\langle \beta_S ,1\rangle =0$ for nearly all
other $S$ by using spectral sequences to prove that the trivial-isotypic
piece of $H(\Pi_n^S)$ vanishes.  The middle ground that is not covered
by our results seems fairly subtle.  Section 4 generalizes the partitioning
of $\Delta (\Pi_n^*)/S_n$ to $\Delta (\Pi_n^*)/S_{\lambda }$ in order to 
prove a second conjecture of Sundaram, regarding representations of
$S_1\times S_{n-1}$ on homology.

\section{Partitioning results}
In Sections 2 and 4, we
will study the flag $h$-vector for $\Delta (\Pi_n )/S_n$ and 
$\Delta (\Pi_n)/S_{\lambda }$, respectively, using
partitionings which express flag $h$-vector coordinates in terms of ascents
and descents in a chain-labeling for orbits of saturated chains in the 
dual poset; this virtually 
necessitates the use of ranks in the dual poset within 
these proofs, despite the fact that
related results and conjectures of Sundaram and our own
spectral sequence arguments are phrased in terms of the ranks
of $\Pi_n$ instead of $\Pi_n^*$.  In an effort to minimize confusion in
converting back and forth between the rank sets for the 
partition lattice and its dual, we will denote by $S^*$ the 
rank set in $\Pi_n^*$ which translates to rank set $S$ in $\Pi_n$ (or
equivalently to corank set $S^*$ in $\Pi_n$).  For the sake of consistency,
the statements of all results will be in terms of rank sets $S$; however, all 
of the proofs in Sections 2 and 4 (as well as Theorem 3.1) work 
internally with 
ascents and descents in the dual poset, so we systematically refer to rank 
sets $S^*$ inside these proofs.  Our arguments also may sometimes
abuse notation by referring to chains when we always mean orbits
of chains.

\begin{figure}[h]
\begin{picture}(150,60)(0,0)
\psfig{figure=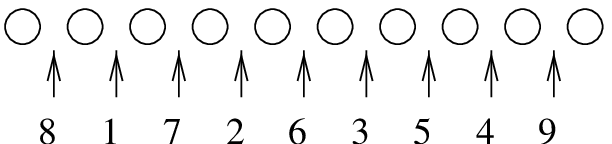,height=1.5cm,width=6cm}
\end{picture}
\caption{Facet achieving $D^3 D^4 A$}
\label{sheila_even}
\end{figure}

Let us depict the facets in $\Delta(\Pi_n^*)/S_n$ (namely the $S_n$-orbits of 
saturated chains in $\Pi_n^*$) by diagrams consisting of 
$n$ balls with bars separating them (or arrows indicating where
these bars are to be inserted) and the numbers from 1 to $n-1$ labeling
the bars (or the arrows).
The bar labels indicate the ranks in $\Pi_n^*$ (or equivalently the
coranks in $\Pi_n $) at 
which the bars are inserted in the course of progressively
refining a single block of $n$ objects into $n$ singleton blocks.
The balls represent
the numbers $1,\dots ,n$ being partitioned, since we may freely permute these
$n$ numbers without switching orbit.  
Figure ~\ref{sheila_even} gives an example (to be used again later)
which begins by refining a block
of size 10 into children of sizes 2,8 and next refines the block of size 8
into children of sizes 2,6.

Each refinement step takes an ordered partition (with the block ordering
coming inductively from the choice of root for the chain orbit) and 
splits one of its blocks into two smaller blocks by inserting a bar
into the block.  Thus, we preserve the order of the original blocks and must
only choose which of the two new blocks goes to the left of the other in the
former position of the parent block.  We make 
the convention of placing each bar as far to the left as possible among all
choices that would give the same saturated chain orbit.
In particular, placing a bar
at position $i$ in a block of size $n$ is 
equivalent to placing the bar at position $n-i-1$ in that block, and we
choose the position farther to the left.  The other situation in which 
there is a choice to make is when there are equivalent blocks of the type
to be split, in which case we refine the leftmost such block.  This only
happens if the blocks 
have the same size and were created from the same parent
at the same step.  

The partitioning in
[He] uses a chain-labeling in which each covering relation $u\prec v$
is labeled by a triple $(i,w,r)$
consisting of the position $i$ of the bar being inserted, 
the word $w$ recording the positions of all bars in $v$
and the rank $r$ at which the block being 
split at step $u\prec v$ was itself
created.  In [He], there is also some sorting of equivalent blocks just prior
to the labeling of each covering relation,
but we may safely ignore this
because our arguments based on the partitioning will only consider 
facets in which no sorting takes place (or in a few cases ranks at which
sorting does not occur 
within facets that do require sorting elsewhere).
We say that a saturated chain orbit in $\Pi_n^*$ has a {\bf topological descent} 
at rank $i$ if the pair of 
covering relations $u\prec v$ and $v\prec w$ labeled $(i,w,r)$ and 
$(i',w',r')$, respectively, satisfy
any of the 
following conditions:
\begin{enumerate}
\item
the bar inserted by $u\prec v$ is farther to the right than the bar insertion
from $v\prec w$ (in which case $i>i'$, so the labels decrease);
\item
the bar insertions $u\prec v$ and $v\prec w$ proceed from left 
to right splitting
a single block of $u$ into three smaller blocks
with the left child from the $u\prec v$ refinement 
strictly larger than the left child from the $v\prec w$ refinement step;
\item
the bar insertions $u\prec v$ and $v\prec w$ refine a single block of $u$ into
three children such that the left children resulting from the $u\prec v$
and $v\prec w$ covering relations both have
size two and the latter gets refined to singletons before the former.
\end{enumerate}

All other ranks in the saturated chain orbit are called (topological)
ascents.  When the above chain-labeling is used to lexicographically order 
facets, 
the topological descents are the ranks which may be omitted from a facet 
to obtain 
codimension one faces that also belong to lexicographically
earlier facets.  For most facets $F_i$, the ranks of the topological
descents in $\Pi_n^*$ are exactly the ranks included in the minimal face
$G_i$ in the 
partitioning for $\Delta(\Pi_n^*)/S_n$, so our aim will be to 
find facets achieving various topological descent sets.

To be more precise, the support of $G_i$ is exactly the ranks of the 
topological descents in $F_i$ if $F_i$ satisfies the nontrivial, non-equal
block condition, as stated just prior to Theorem 3.1.
We should remark that facets $F_i$ in $\Delta (\Pi_n^*)/S_n$ 
violating the nontrivial, non-equal block condition 
have minimal faces $G_i$ whose support is not exactly the set of topological 
descents in $F_i$, but our $b_S(n)>0$ results only need to make use
of facets which do satisfy the non-equal block condition;
when we give partitioning proofs of $b_S(n)=0$ and 
$b_S(n)=b_S(m)$ results, we must consider all facets, but then we use the 
fact that even for facets $F_j$ violating the non-equal block condition, 
the ranks not involving equal blocks are in $G_j$ if and only if they are
topological descent ranks.
These ranks which are governed by descents are enough to show that  
these facets do not have minimal faces of support $S$ forbidden in the
$b_S(n)=0$ results and these ranks also suffice to set up the 
bijection needed in the
new proof of Stanley's $b_S(n)=b_S(m)$ result.

Let us represent the 
rank set $S^* =\{ i_1,i_2,\dots ,i_r\} $ by the word 
$w(S^*)\in \{ A,D\}^{n-2}$ which has a $D$ (for descent) at each position
$i\in S^*$ and an $A$ (for ascent) at each of the remaining position. 
This reflects the fact that in a lexicographic shelling
a saturated chain having descents at exactly the positions in $S^*$ 
would increase the value of the flag $h$-vector coordinate $h_S(\Delta 
(P))$ by one.  Now let us turn to
the following two conjectures of Sundaram [Su]:

\begin{enumerate}
\item
If $1\not\in S $, then $b_S(n) \ne 0$.
\item
Let $b'_S(n)$ be the multiplicity of the trivial representation in
the representation of $S_{n-1}\times S_1$ on homology of the rank-selected
partition lattice $\Pi_n^S$.
If $S=\{ 1,\dots ,i \} $ then $b'_S(n)=1$ and otherwise $b'_S(n)>1$.
\end{enumerate}

\medskip
Sundaram proved in [Su] the first part of her second conjecture, namely that
$b'_S(n)=1$ for $S=\{ 1,\dots ,i\}$.
Theorem 2.1 will confirm the first conjecture.  We
defer the proof of the second conjecture until Section ~\ref{otherlambda}
because it will rely on a partitioning for $\Delta(\Pi_n^*)/S_1\times 
S_{n-1}$ that is given in that section.  Theorem 2.1 is followed by a 
slight generalization, and then we give short new proofs based on 
partitioning for results of Sundaram [Su] and Stanley [St1].

\begin{thm}\label{conj1}
If $1\not\in S$ then $b_S(n)>0$.
\end{thm}

\proof
The requirement $1\not \in S$ translates to $n-2 \not\in S^*$.
Let us exhibit facets in $\Delta (\Pi_n^*)/S_n$
achieving descent set $S^*$ for each such pair $(S,n)$.  First we show how to 
achieve the word $w(S^*) = D^{n-3} A$ for any $n$, handling the cases of 
$n$ even and odd separately.
After this, we will show how to achieve concatenations of such words, 
so as to obtain any word $w(S^*)$ ending in an ascent.  

When $n=2k$, we achieve $D^{k-2}D^{k-1}A$ by first inserting 
$k-1$ bars from left to right sequentially into even positions
$2,4,\dots ,2k-2$, then sequentially inserting bars from right to left 
into odd positions $2k-3,2k-5,2k-7,\dots ,1$; finally, we obtain an
ascent by concluding with a bar in position $2k-1$.  Figure
~\ref{sheila_even} gives an example for $k=5$.
Note that the first $k-2$ pairs of consecutive
bar insertions are topological descents 
because the codimension one face which skips from the bar insertion into 
position $2i$ directly to partition with additional bars at positions
$2i+2,2i+4$ also belongs to the lexicographically earlier facet 
which reverses the order (later in the chain) in 
which bars are inserted into positions 
$2i+1$ and $2i+3$.  The next $k-1$ pairs of consecutive bar insertions 
proceed from right to left and hence are also topological descents.  
The final pair of consecutive bar insertions into positions $1, 2k-1$ gives a
topological ascent. 

\begin{figure}[h]
\begin{picture}(150,45)(0,0)
\psfig{figure=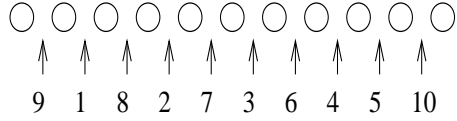,height=1.5cm,width=6cm}
\end{picture}
\caption{Facet achieving $D^4D^4A$}
\label{sheila_odd}
\end{figure}

Similarly, for $n=2k+1$, insert bars sequentially 
into even positions
$2,4,\dots ,2k-2$ then odd positions $2k-1,2k-3,2k-5,\dots ,3,1$ and
finally into even position $2k$, as in Figure ~\ref{sheila_odd}.
The only change is to note that the pair of consecutive insertions at 
positions $2k-2, 2k-1$ is a topological descent since there is a 
lexicographically earlier facet which instead inserts bars first into 
position $2k-3, 2k-1$ and shares a codimension one face skipping rank
$k-1$ with our facet.

To achieve any sequence of ascents and descents that ends in an ascent,
note that the words $D^{n_1}A$ and $D^{n_2}A$ may be concatenated as follows:
\begin{figure}[h]
\begin{picture}(150,60)(0,0)
\psfig{figure=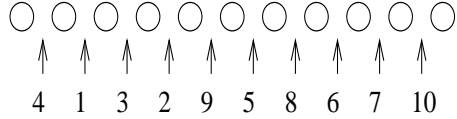,height=1.5cm,width=6cm}
\end{picture}
\caption{Facet achieving $D^3AD^4A$}
\label{sheila_glue}
\end{figure}
apply the above construction for $D^{n_1}A$ 
placing bars into the leftmost $n_1$ available positions, then rather
than concluding with a bar insertion into 
position $n_1+1$, instead place this bar in
position $n_1+2$, creating the leftmost block of size 2 in the beginning of 
the construction for $D^{n_2}A$ using the remaining bar positions.  If there 
are more than two words to concatenate, proceed greedily in this fashion
from left to right among available bar positions.  
Figure ~\ref{sheila_glue} gives an example for $D^3AD^4A$.
\EOP

\medskip 
Now let us strengthen Theorem ~\ref{conj1}.  As before. we state the result
in terms of rank set $S$ in $\Pi_n$, but then use (orbits of) 
saturated chains in $\Pi_n^*$ in the proof.

\begin{thm}\label{conj1+}
Let $S=\{ 1,\dots ,i,j_1,\dots ,j_l\} $ for $j_1-i>1$.  If $i\le l$, then
$b_S(n)>0$.
\end{thm}

\proof
First let us consider the case when $i=|S|/2$.  This means $S^*$ has
exactly $i$ descents interspersed with ascents and then has at
least one ascent immediately preceding the final 
string of $i$ consecutive descents.
To realize this pattern in a saturated chain in $\Pi_n^*$, 
begin by inserting bars from left to right 
creating blocks of size 1 or 2, making a block of size 2 at
each descent and a block of size 1 for each ascent.
This accounts for the first $i$ descents and the ascents in which they
are interspersed.  Now we refine the remaining
$i+ 1$ blocks of size 2 
from right to left to achieve the last $i$ descents.  
Notice that there must be one extra block of size 2 immediately to the 
right of the rightmost block of size one that we have created, because the
size of $n$ dictates that there must be $i+1$ further refinements, forcing
this rightmost block to have size 2.

Now suppose $l>i$.  Let us then break $w(S^*)$ into subwords $w_1,w_2$ such
that $w(S^*) = w_1\circ w_2$ and the word $w_1$ has 
exactly $l-i$ descents including a terminal one.
If $w_2$ begins with an ascent, then we achieve $w_1A$ by the 
construction for words ending in an ascent applied to the leftmost bar
positions, and then we achieve
$w_2'$ such that $w_2 = Aw_2'$ by the construction of the previous 
paragraph on the remaining bar positions.  
On the other hand, if $w_2$ begins with a descent,
then we use the 
previous theorem's construction for the maximal word $w_1'$ such that
$w_1 = w_1' \circ D^{r_1}$, since $w_1'$ must end in an ascent; then we
achieve $D^{r_1 + r_2}Aw_2''$ such that $w_2 = D^{r_2} Aw_2''$ by the 
following construction:
if $r_1 $ is even, then insert bars sequentially from left to 
right into even positions $2,4,\dots ,2r_2 + r_1 $, then from right to 
left into the odd positions $2r_2 + r_1 -1,2r_2 + r_1 -3, \dots ,
2r_2 + r_1 - (r_1-1)$, then use the $i=|S|/2$ procedure for 
the bar positions to the right of $2r_2 + r_1$ and conclude by placing bars
right to left into consecutive odd positions $2r_2 + r_1 - (r_1+1),\dots ,1$;
for $r_1$ odd, sequentially insert bars into even positions 
$2,4,\dots ,2r_2 + r_1-1$, then insert bars right to left in odd 
position $2r_2 + r_1,2r_2 + r_1 -2, \dots ,2r_2 + 1$, and finally use the 
$i=|S|/2$ procedure on the positions to the right of 
$2r_2 + r_1$.
\EOP

\medskip
Sundaram [Su, p. 288] showed that $b_S(n)=0$ whenever any of the 
following conditions are met:\\
\begin{enumerate}
\item
$S=\{ 1,\dots ,i \} $ with $i>0$.
This result of [Ha] and [Su] is recovered by partitioning in [He].
\item
$[1,\lfloor (n+1)/2 \rfloor ] \subseteq S$.
\item
$S=[1,r]\cup a$ for $a\not\in [{r+2\choose 2},n-r-1]$. 
\item
$S=[1,r] - k$ for $n$ even and $k=n/2-1$, provided 
$k=\frac{n}{2} - 1 \le r\le n-4$.  
\end{enumerate}

\medskip
Sundaram asked (private communication) if these results could 
easily be recovered using the
partitioning for $\Delta(\Pi_n)/S_n $.
We now give proofs by partitioning for these 
results and then provide further results about when $b_S(n)>0$ and 
when $b_S(n)=0$.  Let us begin with Item 2.

\begin{thm}
If $[1,\lfloor (n+1)/2 \rfloor ] \subseteq S$, then $b_S(n)=0$.
\end{thm}
\proof
Consider the 
partition immediately before the final string of descents in the orbit of
any saturated chain in $\Pi_n^*$.  At this point,
less than half the bars have been inserted, so some blocks have size 
larger than 2.  Each such block forces an ascent, a contradiction.  
\EOP

\medskip
Next we give a slight strengthening of Sundaram's fourth result.

\begin{thm}
If $S=[1,r] - k$ for $k> r/2$, then $b_S(n)=0$.
\end{thm}
\proof
Suppose the orbit $O(C^*)$
of some saturated chain in $\Pi_n^*$ achieves the set of 
coranks $S^*$ for
$S=[1,r]-k$, i.e. suppose it achieves $w(S^*)=
A^{n-2-r}D^{r-k}AD^{k-1}$.  For $O(C^*)$
to begin with $n-2-4$ ascents, 
bars must be inserted left to right creating blocks of nondecreasing size.  
Each of these blocks of
size larger than 2 will necessitate an ascent to complete its refinement some
time after the first descent.
Thus, there may be at most one block of size larger than 2
created by the initial string of ascents.  In addition, any
blocks of size 2 which are created initially must later be
split from left to right for none of the initial bar insertions
to be topological descents.

Thus, the initial string of ascents creates at 
most one block of size 2, so it creates some number of trivial
blocks, followed by at most one block of size 2, and then at most one
larger block.  The first topological descent must come from proceeding right
to left in order to refine 
the unique block of size two or else we get a topological descent by
splitting off the leftmost singleton in a block of size larger than two.  
In either case, this step must be followed by an 
ascent.  Now we need to achieve a nonempty string of descents to complete the
refinement, but it is impossible to 
completely refine the block which had size larger than two using only 
topological descents, by the same 
argument that was used to show $b_{1,\dots ,i}(n)=0$ in [He].
\EOP

\medskip
Another application of the partitioning is 
a simple proof for the following result of 
Stanley [St1].

\begin{thm}
Let $S=\{ j_1,\dots ,j_l \} $.  Then $b_S(m)=b_S(n)$ if $m,n > 2j_l$.
\end{thm}

\proof
A saturated chain in $\Pi_n^*$ achieving a set $S$ as above 
must begin with $n-j_l-2$
consecutive ascents, so 
we must insert bars left to right creating blocks of nondecreasing size.  
These ascents create at most 
$j_l$ blocks of size greater than one since there are only $j_l$ remaining
refinement steps to completely refine these blocks.  Thus, the initial 
ascents must create only blocks of size one from left to right for the 
first $n-2j_l -2$ bar insertions.
We get a bijection between facets
in $\Delta(\Pi_n^*)/S_n$ and in $\Delta(\Pi_m^*)/S_m$ which contribute 
minimal faces of support $S$ to their respective partitionings, 
by changing the number of initial blocks of size 
one and otherwise letting the bar insertions agree once we have split off
the necessary number of singletons from each facet.
\EOP

\section{Conditions under which $b_S(n)=0$}

In this section, we will recover Sundaram's third result using 
spectral sequences, and then we generalize her result by 
replacing the single rank $a$ by 
a collection of ranks which are disjoint from the consecutive initial 
ranks $1,\dots ,i$.
First we use a partitioning for $\Delta(\Pi_n^*)/S_n$ to obtain conditions
under which $b_S(n)>0$, before using spectral 
sequences to show $b_S(n)=0$ for nearly all other $S$.   

Let $S= \{ 1,\dots ,i,j_1,\dots ,j_l \} $ with $j_1-i>1$, so $S^* = 
\{ n-1-j_l, n-1-j_{l-1},\dots ,n-1-j_1\} \cup [n-1-i,n-2] $.  
Let $E(C^*)$ denote the lexicographically smallest extension
of a face $C^*$ in $\Delta (\Pi_n^*)/S_n$ to a saturated chain orbit, 
based on the chain-labeling of [He].  
Let us say that a 
chain $C=\alpha_1\prec\cdots\prec \alpha_i < \beta_1 < \cdots < 
\beta_l $ of support $\{ 1,\dots ,i,j_1,\dots ,j_l \} $ satisfies the 
{\bf non-equal block condition} if the extension $E(C^*)$ to a saturated
chain does not have any pairs of equal blocks created from
the same parent either in a single refinement step or in consecutive 
refinement steps.  If we relax this requirement to allow equal blocks of
size two, we call this the {\bf nontrivial, non-equal block condition}.
Theorem 3.1 is not tight in
that there are rank sets $S=\{ 1,\dots ,i,j_1,\dots ,j_l \} $ 
such that $b_S(n)>0$
but where every facet whose minimal face has support $S$ violates the 
non-equal block condition; the
situation seems much more subtle when one removes the non-equal block
condition.  Let $Stab(C)$ denote the stabilizer of a chain $C$.

\begin{thm}
Let $C=\alpha_1\prec\cdots\prec\alpha_i < \beta_1 < \cdots < \beta_l$ be a
chain in $\Pi_n $ of support $S=\{ 1,\dots ,i,j_1 , \dots , j_l \} $ for
$j_1 > i+1$ that satisfies the non-equal block condition.  
Furthermore, suppose that $\alpha_i$ has $i$ nontrivial blocks 
$B_1,\dots ,B_i$ of size 2 and that $\beta_1 $ has $i+1$ nontrivial blocks
$C_1,\dots ,C_{i+1}$ belonging to distinct $Stab(\beta_1 < \cdots < 
\beta_l)$-orbits such that $B_r \subseteq C_r $ for $1\le r\le i$.
Then $b_S(n)>0$.
\end{thm}

\proof
Let $C$ be a chain as above.  Either
$E(C^*)$ will have (topological) descents at exactly the ranks in $S^*$ or
we will construct a closely related chain $C'$ with the desired topological
descent set for $E(C')^*$ such that $C'$ satisfies the 
non-equal block condition.  Thus, either $E(C^*)$ or $E((C')^*)$ will 
contribute to the partitioning a minimal face of 
support $S^*$, implying
$b_S(n)>0$.

The saturated chain $E(C^*)$ is obtained 
by extending each interval $u<v$ in $C^*$
by inserting bars left to right and splitting each block of $u$ from left
to right into nondecreasing pieces, so that each rank in the extension is
a topological ascent.
Thus, the topological descents of $E(C^*)$ are a subset of the ranks
in $C^*$, i.e. they are a subset of $S^*$.  Furthermore, $E(C^*)$ must have
topological descents at the topmost $i$ ranks because $C^*$ may be chosen to
conclude by refining from right to left the $i+1$ 
blocks of size 2 which are children of the $i+1$ 
blocks in $\beta_1$ that are in distinct $Stab(\beta_1 < \cdots < \beta_l)$
orbits.  As an example, the chain orbit $C^* = (00|000|00000 
< 00|0|00|0|0|0|0|0 < 00|0| 0|0 |0|0|0|0|0 < \hat{1} )$ 
of support $2,7,8$ extends to a chain $E(C^*)$
which sequentially inserts bars in positions $2,5,3,6,7,8,9,4,1$, and this
has descents at ranks $2,7,8$; it contributes to flag $h$-vector
coordinate $h_{2,7,8}$ for $\Delta(\Pi_{10}^*)/S_{10}$ or equivalently to 
$h_{9-2,9-7,9-8} = h_{1,2,7}$ for $\Delta(\Pi_{10})/S_{10}$. 

Let us next consider the lower ranks in $E(C^*)$ and specifically
how to turn ascents at ranks in $C^*$ into topological descents.  
If there is an ascent at a rank $u\in E(C^*)$, then the bar
inserted at the covering relation $t\prec u$ in $E(C^*)$
is to the left of the bar inserted at the covering relation $u\prec v$ 
in $E(C^*)$, and the two bars are either (1) 
inserted in different 
blocks or (2) inserted into a single block creating children 
$B,B',B''$ from left to right such $|B| \le |B'|$.  

For the moment, let
us assume that $t$ and $v$ are 
not in $C^*$, so that our goal will be to replace the ascent at $u$ by a 
descent while preserving the ascents $t$ and $v$.
For ascents of type (1), we obtain from $C^*$ a new chain $(C')^*$ 
by inserting the right bar before the left one.  Note that the necessary
ascents are preserved
since the label leading upward to $u$ was increased in value and
the label leading upward from $u$ was decreased;
furthermore, nonequivalent blocks in $C^*$ are still nonequivalent 
in $(C')^*$, preserving the requirement about $Stab (\beta_1 < \cdots < 
\beta_l)$ orbits.  For ascents of type (2), we replace the consecutive
bar insertions which give left children $B, B'$ such that $|B|\le |B'|$ 
by left to right bar insertions instead sequentially yielding left children 
$B',B$ to produce a topological descent.  
In the remainder of $(C')^*$, we refine $B,B'$ just as $C^*$ 
would, though this could in theory impact later ascents and descents at
ranks in $C^*$ since now $B'$ is to the left of $B$.  We can choose our 
chain to avoid turning ascents to descents and vice-versa
as long as we proceed from lower to higher ranks in creating $C'$.
This modification of $C^*$ into $(C')^*$ for type 2 
also gives ascents immediately above and below the descent at $u$,
because in order 
for the bar creating $B$ to be farther to the right than the bar inserted 
just before it,
the bar which instead creates $B'$ must also be farther to the right, and 
since $|B'|> |B|$, we observe that $B'$ is also larger than any block
$B''$ created just prior to $B$ from the same parent as $B$, since we would
have $|B''| < |B|$.  Similarly, we are 
assured of an ascent immediately after the descent at $u$, and again we
have preserved block-nonequivalence as needed.
Also, the left child is never larger than the right 
child in a refinement since $|B'|\le |B''|$ and $|B|\le |B''|$.

Now let us more generally consider the possibility 
that the ascent $u$ is among $r$ consecutive
elements $u_1,\dots ,u_r \in C^*$ for $r \ge 1$.  Let us describe how to
obtain $(C')^*$ which has topological descents at all of 
these $r$ consecutive 
ranks and ascents immediately above and below them.  Let $u_0\prec u_1$
and $u_r\prec u_{r+1}$ be the covering relations of $E(C^*)$ immediately
below and above the $r$ consecutive ranks.
In $C'$, we refine the set of blocks in $u_0$ from right to left
and within each block of $u_0$
insert bars left to right with new blocks decreasing
in size from left to right.  Thus we get a string of descents.  This
is immediately preceded and followed by ascents, since (just as in the 
$r=1$ case) the label for $u_1\prec u_2'$ in $E(C')^*$ is smaller than 
the corresponding label in $E(C^*)$, and the label for $u_{r-1}'\prec u_r$
in $E(C')^*$ is larger than the corresponding label in $E(C^*)$.  Our
conversion of $C$ to $C'$ also preserves block-nonequivalence as needed.
\EOP

\medskip
Next we use a spectral sequence for a filtered complex to give
upper bounds on $b_S(n) = \langle 1,\beta_S(\Pi_n)\rangle $ by showing 
that the trivial-isotypic piece of  
$E^2(\Delta(\Pi_n^*))$ vanishes 
except when a certain type of chain of support $S$ exists.  See [We]
for background on spectral sequences and in particular on how they give
(upper) approximations on homology.  We begin with
a new proof that $b_{\{ 1,\dots ,i\} }(n)$ 
which we then generalize from rank set $[1,r]\cup \{ a\} $ 
to rank set $[1,r]\cup \{ a_1\cup \cdots\cup a_l \} $ for $l\ge 1$.
Now we come to Sundaram's third result.

\begin{thm}\label{onelone}
Let $S=[1,r]\cup a$ for 
$a\not\in [{r+2\choose 2},n-r-1]$.  Then $b_S(n)=0$.
\end{thm}
\proof
Consider the rank selection $\{ 1,2,\dots ,i,j\} $ for any 
$i < j < {i+2 \choose
2}$.  Let $C$ be the set of chains in $\Pi_n $ which are supported on a 
subset of these ranks.

For each chain $\Gamma $ in $C$, define $f(\Gamma )$ to be 
$2 J(\Gamma ) + I(\Gamma )$ where 
$J(\Gamma ) = 1$ if $\Gamma $ contains an 
element of rank $j$, and $J(\Gamma ) = 0$ otherwise; $I(\Gamma )$ is 
defined to be the rank of the maximal element of $\Gamma $ which has rank
at most $i$.

Note that $\bdry (\Gamma )$ is a linear combination of chains $\Gamma '$
with $f (\Gamma ') \le f(\Gamma )$.  So $f$ is a filtering on the 
complex $(C,\bdry )$ and so we can approximate the homology $H(C,\bdry )$
with $E^1 = H(C,\bdry^0 )$ where $\bdry^0 $ is the piece of the boundary
$\bdry $ which is fixed by the filtration function.  

It is easy to see what $\bdry^0$ does: if $\Gamma $ has the form 
$$\Gamma = \hat{0} < x_1 < \cdots < x_s < y < \hat{1}$$
where $y$ is of rank $j$, then
$$\bdry^0 (\Gamma ) = \sum_{l=1}^{s-1} (-1)^{l-1} \bigg\{\hat{0} < \cdots < 
\hat{x_l} < \cdots < x_s < y < \hat{1}\bigg\} .$$
If $\Gamma $ is of the form $\hat{0} < x_1 < \cdots < x_s < \hat{1}$ where
$\rm{rk} (x_s ) \le i$, then 
$$\bdry^0 (\Gamma ) = \sum_{l=1}^{s-1} (-1)^{l-1} \bigg\{\hat{0} < \cdots < 
\hat{x_l} < \cdots < x_s < \hat{1}\bigg\} .$$

By examining the form of $\bdry^0$, one sees that $(C,\bdry^0)$ can be split as a direct sum
$$ (C,\bdry^0 ) = \plus_{\alpha\in R^{\le i}} C([\hat{0},\alpha ],\delta )
\oplus \plus_{\gamma\in R^{\le i}\atop \beta\in R^j , \gamma < \beta }
C([\hat{0},\gamma ],\delta ).\eqno(3.1) $$
Here $C([\hat{0},\alpha ],\delta )$ denotes the usual order complex of 
the poset $[\hat{0},\alpha ]$.  Let $R^{\le i}$ denote the set of poset 
elements of rank at most $i$, and let $R^j$ denote the set of poset elements
of rank exactly $j$.

Note that $f(\sigma\Gamma ) = f(\Gamma )$ for $\sigma \in S_n$.  Therefore
$S_n$ commutes with the boundary $\bdry^0$ and so it makes sense to talk
about the $S_n$-module structure of the complex $(C,\bdry^0)$.  The 
$S_n$-module structure is best described in pieces corresponding to the 
two major summands in 3.1.

The summand $\plus_{\alpha\in R^{\le i} } C([\hat{0},\alpha ],\delta )$
corresponds to the space of all chains which have no element of rank $j$.
Likewise, the summand $\plus_{\alpha\in R^{\le i}\atop \beta\in
R^j, \alpha < \beta } C([\hat{0},\alpha],\delta )$ corresponds to 
the span of all chains which do have an element of rank $j$.  These two
subspaces are $S_n$-invariant.

For the first summand, let $I = \{\alpha_1 ,\dots ,\alpha_l \} $ be the
set of representatives from the orbits of $S_n$ acting on $R^{\le i}$.
Then as an $S_n$-module, the first summand is
$$  \plus_{\alpha \in I} ind_{Stab (\alpha )}^{S_n} (
C([\hat{0},\alpha ] , \delta )) \eqno(3.2)$$
where $Stab(\alpha )$ denotes the stabilizer of $\alpha $ in $S_n$.

For the second summand we have a similar description.  Let $J = \{
\gamma_1 < \beta_1, \gamma_2 < \beta_2 , \dots ,\gamma_m < \beta_m \} $
be a set of representatives from the orbits of $S_n$ acting on the 
set $\{ \gamma < \beta : \rm{rk} (\gamma ) \le i , \rm{rk} (\beta ) = j, 
\gamma < \beta \} $.  Then the second summand is 
$$  \plus_{ \{ \gamma < \beta \} \in J } ind_{Stab(\gamma < 
\beta )}^{S_n} ( C([\hat{0},\gamma ],\delta )).\eqno(3.3) $$

So, we can write the following expression for $E^1 = H(C,\bdry^0 )$.
$$  E^1 = \plus_{\alpha\in I } ind_{Stab(\alpha )}^{S_n} 
(H([\hat{0},\alpha ]))
\oplus \plus_{ ( \gamma < \beta ) \in J} ind_{Stab(\gamma < \beta )}^{S_n}
(H([\hat{0},\gamma ])).\eqno(3.4) $$

Our goal is to show that the multiplicity of the trivial representation in
$E^{\infty }$ is 0.  We will begin by characterizing the trivial-isotypic
component in $E^1$.

As a notational convention, whenever $V$ is a $G$-module for any group $G$, 
let $V^G$ denote the trivial-isotypic component of $V$.  By Frobenius
Reciprocity,
$$  (E^1)^{S_n} \approx \plus_{\alpha\in I} 
H([\hat{0},\alpha] )^{Stab(\alpha )} \oplus \plus_{(\gamma < \beta )\in J}
H([\hat{0},\gamma ])^{Stab(\gamma < \beta )}.\eqno(3.5) $$

Let $\alpha = A_1 | \cdots | A_k | B_1 | \cdots | B_l|C_1|\cdots |C_m$
where the $A_u$ all have size 1, the $B_v$ all have size 2 and the $C_w$ 
all have size greater than 2.  Corresponding to this decomposition, 
$$  St(\alpha ) = \prod_{u\ge 1} (S_{m_u} \wr S_{u })\eqno(3.6) $$
where $m_u$ is the number of blocks of $\alpha $ of size $u$ and $S_{m_u}
\wr S_u$ denotes a wreath product of symmetric groups.  Likewise,
$$  H([\hat{0},\alpha ]) \cong  \tensor_{u\ge 3 }
H (\Pi_u )^{\otimes m_u}.\eqno(3.7) $$

The action of $Stab(\alpha )$ on $H([\hat{0},\alpha ])$ is given by
an action of each $S_{m_u} \wr S_u $ on the tensor factor 
$H(\Pi_u )^{\tensor m_u} $.  This wreath product action is the one in which 
the $m_u$ copies of $S_u$ act on $H(\Pi_u )$ in the usual way.  The overlying 
copy of $S_{m_u}$ acts according to the trivial representation if $u$ is odd
and the sign representation if $u$ is even.  

From this description of the action of $Stab(\alpha )$ on 
$H([\hat{0},\alpha ])$, we will deduce that 
$$  H([\hat{0}, \alpha ])^{Stab(\alpha )} = \complex \eqno(3.8) $$ 
if $\alpha $ has a single 
block of size 2 and all other blocks of size 1, and that
$H[\hat{0},\alpha ])^{Stab(\alpha )}$ is 0 otherwise.  
To obtain 3.8, we use the well-known fact that
$$ H(\Pi_u)^{S_u} =  \bigg\{ { 0 \hspace{.1in}\rm{ for }
\hspace{.1in}u\ge 3 \atop \hspace{.11in} \complex 
\hspace{.1in}\rm{ for } \hspace{.1in} u=1,2} \eqno(3.9) $$
Notice that the 
trivial representation of $Stab (\alpha )$ is $\tensor_{u\ge 1} 
1_{m_u} \wr 1_u$, where $1_u $ denotes the trivial representation of
$S_u$.  Computing inner products, we see that a $Stab(\alpha )$-representation
$\chi \wr \psi $ will
not contain the trivial representation unless both $\chi $ and $\psi $
do as well (cf. [JK, chapter 4]).  Thus, 3.8 follows from 3.9 along with 
our above description of the action of $Stab (\alpha )$ on 
$H([\hat{0},\alpha ])$.

Similar reasoning allows us to analyze the summand
$\plus_{(\gamma ,\beta)\in J} H([\hat{0},\gamma ])^{Stab(\gamma < \beta )}$
from (3.5).
However, there is a subtlety here in that $Stab(\gamma < \beta )$ is not 
necessarily the full automorphism group of $\gamma $.  For each block $B$ of
$\gamma $, the automorphism group $Stab(\gamma < \beta )$ will certainly 
contain $S_B$.  However, different blocks of $\gamma $ of the same size may
not be interchanged by $Stab(\gamma < \beta )$ because they reside in blocks
of $\beta $ which have different size.  The conclusion is that
$H([\hat{0},\gamma ])^{Stab(\gamma < \beta )} = 0$ unless every nontrivial 
block of $\gamma $ has size 2 and if $U,V$ are blocks of $\gamma $ having
size 2, then $U$ and $V$ are contained in blocks of $\beta $ which have 
different sizes.

Let $\gamma < \beta $ be a pair for which $H([\hat{0},\gamma ])^{Stab (\gamma
< \beta )} $ is nonzero.  We will identify the structure of $H([\hat{0},
\gamma ])^{Stab(\gamma < \beta )} $ more explicitly.  Let $U_1,\dots ,U_t$
be the nontrivial blocks of $\gamma $ (all of which have size 2).  For each
$l$, let $V_l$ be the block of $\beta $ which contains $U_l$.  We know that
$|V_l|=|V_m|$ implies $l=m$.

The poset $[\hat{0},\gamma ]$ is isomorphic to the Boolean algebra $B_t$ and 
so we know that the module
$H([\hat{0},\gamma ])$ has dimension 1 in degree $t$ and 
dimension 0 in all other degrees.  In addition we can explicitly give the 
homology representative $\rho [U_1,\dots ,U_t]$ in degree $t$:
$$  \rho [U_1,\dots ,U_t] = \sum_{\sigma\in S_t} sgn (\sigma )
\{ \hat{0} < U_{\sigma 1} < U_{\sigma 2} | U_{\sigma 3} < \cdots < 
U_{\sigma 1} | U_{\sigma 2} | \cdots | U_{\sigma (t-1) } < \gamma \} .
\eqno(3.12) $$

Consider now the next step in the spectral sequence.  We are going to compute
$E^2 = H(E^1,\bdry^1)$ where $\bdry^1$ is the differential induced on $E^1$
by the piece of the original boundary which reduces the filtration 
function by 1.  By the definition of $f$, $f(\Gamma') = f(\Gamma ) -1$ for
chains $\Gamma ' \subseteq \Gamma $ iff $\Gamma '$ is obtained from 
$\Gamma $ by removing the maximal element of $\Gamma $ whose rank is in
$\{ 1,\dots ,i\} $.

Referring to (3.8) we see that the $S_n$-invariants in the first summand
of (3.5) consists of the 
single vector
$$  \plus_{type(\alpha )= 2,1^{n-2}} H([\hat{0},\alpha ]) = 
\langle \sum_{type(\alpha )= 2,1^{n-2}} (\hat{0} < \alpha < \hat{1} ) 
\rangle_{ \complex }.\eqno(3.13) $$

Applying $\bdry^1$ to (3.13) 
gives a non-zero multiple of $(\hat{0} < \hat{1})$,
hence the first summand contributes nothing to the kernel of $\bdry^1$ so 
nothing to $E^2$.

Moving to the second summation, let $\gamma < \beta $ be a pair for 
which $H([\hat{0},\gamma ])^{Stab(\gamma < \beta )}$ is nonzero.
As with earlier notation, let $\gamma $ have nontrivial blocks 
$U_1,\dots ,U_t$.  It is straightforward to deduce from equation (3.12) that
$$  \bdry^1 \rho [U_1,\dots ,U_t ] = \sum_{i=1}^t (-1)^i
\rho [U_1,\dots ,\hat{U_i} , \dots , U_t ].\eqno(3.14) $$

Let $V_1^{(1)},V_2^{(1)},\dots ,V_{m_1}^{(1)},V_1^{(2)},
V_2^{(2)},\dots ,V_{m_2}^{(2)},
\dots ,V_1^{(s)},\dots ,V_{m_s}^{(s)}$ be the nontrivial blocks of $\beta $
indexed so that $|V_i^{(j)}| = v_j$ for all $i,j$ and with 
$2 \le v_1 < v_2 < \cdots < v_s$.  Note that $\rm{rk}(\beta ) = \sum_{l=1}^s 
m_l (v_l -1)$.  Since $\rm{rk} (\beta ) < {i+1 \choose 2} = 1 + 2 + \cdots + 
(i+1)$, it follows that $s\le i$.  So there are at most $i$ different 
non-trivial block sizes in $\beta $.

We now compute the contribution to $E^2 = H(E^1 ,\bdry^1 ) $ made by the 
second summation.  By (3.14), for each $\beta $ at rank $j$, $\bdry^1 $
preserves $\plus_{\gamma < \beta } H([\hat{0},\gamma ])^{Stab (\gamma < 
\beta )}.$

Let $z_l = {v_l \choose 2} m_l$ so that $z_l$ is the number of pairs of 
numbers which occur in the same block of size $v_l $ in $\beta $.  Let 
$Z = Z_1 \cup Z_2 \cup \cdots \cup Z_s $ be a vertex set where $Z_l$ 
contains $z_l$ nodes $Z_l = \{ x_1^{(l)},\dots ,x_{z_l}^{(l)}\} .$

Define a simplicial complex $\Delta_{\beta }$ with vertex set $Z$ by 
saying that $\Delta_{\beta }$ contains all subsets $S=\{ s_1, \dots ,
s_t \} $ such that $| S\cap Z_l | \le 1$ for all $l$.

From the description of $H([\hat{0},\gamma ])^{Stab(\gamma < \beta )}$ in 
(3.12) and the formula for the boundary $\bdry^1 $ it is clear that 
$H(\plus_{\gamma } H([\hat{0},\gamma])^{Stab (\gamma < \beta )}, \bdry^1 )
= H (\Delta_{\beta })$ where $H(\Delta )$ is the ordinary simplicial 
homology of $\Delta $.

Using arguments similar to previous ones, 
$$(E^2)^{S_n} = \plus_{\beta\in L} H(\plus_{\gamma } H([\hat{0},\gamma 
])^{Stab (\gamma < \beta )} ,\bdry^1 )^{Stab (\beta )} = \plus_{\beta\in L}
H (\Delta_{\beta } )^{Sym (Z_1 ) \times \cdots \times Sym (Z_s)}$$
where $L$ is a complete set of representatives for the orbits of $S_n$ on
rank $j$ and where $Z_1,\dots ,Z_s$ in the summand $\beta $ are as above.  
But clearly $H(\Delta_{\beta })^{Sym (Z_1) \times \cdots \times Sym (Z_s)}
= 0$ as the projection by the trivial character of $Sym (Z_1)\times \cdots
\times Sym (Z_s )$ maps $\Delta_{\beta }$ to the simplicial complex of all
subsets of $\{ 1,2,\dots ,s\} $ which is acyclic, because $s\le i$.

A slight modification of this argument shows that the rank selected 
homology is 0 for ranks $1,2,\dots ,i, j$ when $j\ge n-i$.  The idea 
is that $\beta $ has at most $n-j \le i$ distinct blocks since there are only 
$n-j-1$ covering relations 
from $\beta $ to $\hat{1} $ in which to merge the blocks of 
$\beta $.
\EOP

\medskip
Next we will extend the above argument, beginning with 
the choice of filter.  If we allow $\Gamma = \alpha_1\prec\alpha_2\prec
\cdots\prec\alpha_k < \beta_1 < \cdots < \beta_l$ above
the $\alpha $ chain, then we let $J(\Gamma )=2l$ and use the appropriately
adjusted stabilizers and $\Delta_{\beta_1 < \cdots < \beta_l}$.

\begin{thm}
If $b_S(n)>0$ for $S=[1,i]\cup \{ j_1,j_2,\dots , j_l \} $ with 
$j_1 > i+1$, then
there exists a chain $\alpha_1\prec\cdots\prec \alpha_i < \beta_1 < 
\cdots < \beta_l$ of support $S$ such that (1) $\alpha_i$ consists of $i$
blocks of size 2, and (2) $\beta_1$ includes nontrivial
blocks $B_1,\dots ,B_{i+1}$ all belonging to 
distinct $Stab(\beta_1 < \cdots < \beta_l)$-orbits.
\end{thm}

\proof
We will mimic the reasoning given in the case 
$\ell = 1$.  First, we define the filtering
$f$.  Let  
$$\Gamma = \alpha_1 < \alpha_2 < \cdots < \alpha_s 
< \beta_1 < \cdots < \beta_t$$
be a chain in the rank selection where the ranks of the 
$\alpha_u$ are in $[1,i]$ and the ranks of the $\beta_v$
are in $\underline{j} = 
\{j_1, j_2, \ldots ,j_{\ell}\}$.  Define $f(\Gamma)
= 2t + \rm{rk}(\alpha_s)$.
 
As in the proof of Theorem 3.2, the $E^1$ term in the spectral
sequence corresponding to this filtration is a direct sum, over
$S_n$-orbits of elements $\alpha_s$ with ranks in $[1,i]$,
and over chains $\beta_1 < \cdots < \beta_t$ with ranks contained
in the set $\underline{j}$ which also satisfy $\alpha_s < \beta_1$.
The summands take the form
$$ind_{Stab(\alpha_s,\underline{\beta})}^{S_n} (H(\alpha_s) )
 \eqno(3.15)$$Here
$Stab(\alpha_s,\underline(\beta))$ denotes the stabilizer of the 
chain $\alpha_s < \beta_1 < \cdots < \beta_t$.
 
We next compute the multiplicity of the trivial character in each
of the summands in $(3.15)$.  As before, we invoke Frobenius reciprocity
to argue that the multiplicity of the trivial character in $(3.15)$
is equal to the multiplicty of the trivial character of 
$Stab(\alpha_s,\underline{\beta})$ in $(3.15)$.  
 
By reasoning similar to that used in 
the proof of Theorem 3.2, the trivial character has multiplicity $0$
in $(3.15)$ unless every 
nontrivial block of $\alpha_s$ has cardinality two.  
Assume therefore that every nontrivial block of $\alpha_s$ has 
cardinality two.  If any pair of these nontrivial blocks are in 
the same orbit of $Stab(\alpha_s,\underline{\beta})$ then again
the multiplicity of the trivial representation is $0$.  

Now consider the summand given by a chain $\underline{\beta} 
= \beta_1 < \cdots < \beta_t$ in
the trivial-isotypic piece of the $E^2$ term of the spectral sequence.
Similarly to when $l=1$, this summand is 
isomorphic to the simplicial homology of a simplicial complex
$\Delta_{\underline{beta} }$ whose ground set consists of the 
$Stab (\beta_1 < \cdots < \beta_t)$-orbits of the size 2 blocks
that are contained in blocks of $\beta_1$.  Let $j$ be the size of
this ground set.  The faces of $\Delta_{\underline{\beta }}$
are the $Stab (\beta_1 < \cdots < \beta_1)$-orbits of chain elements
$\alpha < \beta_1$ of rank at most $i$ whose nontrivial blocks all have 
size 2.  Thus, we are taking the simplicial homology of the $(i-1)$-skeleton
of a $(j-1)$-simplex, so this is 0 unless $j>i$, as desired.  Hence, 
$(E^2)^{S_n}$ vanishes unless there is a chain $\underline{\beta }$ with
$j>i$.  This implies $b_S(n)=0$ unless there is such a chain, as desired.
\EOP

\section{Partitioning $\Delta(\Pi_n)/S_{\lambda }$ and a conjecture of 
Sundaram}\label{otherlambda}

Next we prove the second conjecture of Sundaram [Su], by first giving
a partitioning for $\Delta(\Pi_n)/S_1 \times S_{n-1}$, and more
generally for $\Delta(\Pi_n)/S_{\lambda }$ for any Young subgroup of $S_n$.  
Instead of using bars to partition $n$ balls, now we partition the 
multiset $\{ 1^{\lambda_1},\dots ,k^{\lambda_k} \}$.
Recall that [He] always chose the 
leftmost of equivalent positions in which to insert bars, splitting
a block by inserting a bar with the smaller resulting block to its left;
we more generally need an ordering 
on blocks which are subsets of $\{1^{\lambda_1 },\dots ,k^{\lambda_k }\} $
to decide which block goes to the left of each bar insertion.  The entire
partitioning argument of [He] will go through directly if we use any block
ordering that satisfies the lengthening condition, defined as follows:

\begin{defn}
A block order satisfies the {\bf lengthening condition (LC)} if
$$B \le B' \Rightarrow B \le B B',$$ where $B B'$ denotes the concatenation
of the two blocks, so the multiplicity in $B B'$ of any letter appearing in 
both $B$ and $B'$ is the sum of the multiplicities.
\end{defn}

\medskip
Denote by $w_B$ the word obtained from a block $B$ by rearranging the 
letters in $B$ into increasing order.  If one views blocks as monomials, then
any monomial term order will satisfy the lengthening condition.  However,
the distinguished block order (described below) satisfies the lengthening
condition but is not a monomial term order; we will use the length-lex
order to partition $\Delta (\Pi_n)/S_{\lambda }$ for an arbitrary 
$\lambda $. 

\begin{defn}
In the {\bf length-lex} block order, a block $B_1$ is smaller than
a block $B_2$ if $|B_1| < |B_2|$ or if $|B_1| = |B_2|$ and $w_{B_1}$
is lexicographically smaller than $w_{B_2}$.
\end{defn}

When $\lambda_k=1$, the distinguished block order will give 
a different partitioning that is more convenient for counting 
minimal faces of particular supports in the partitioning and in particular
for proving a second conjecture of Sundaram.

\begin{defn}
Suppose $\lambda $ is a partition in which $\lambda_k =1$, and let 
$s$ be a letter appearing with multiplicity one.
In the {\bf distinguished block order} for 
$\Delta(\Pi_n)/S_{\lambda }$, a pair of blocks $B_1,B_2$
satisfy $B_1< B_2$ if $s\in B_1, s\not\in B_2$ or if $s\not\in B_1, B_2$
and $B_1<B_2$ in the length-lex block order.  The letter $s$ cannot 
belong to two different blocks, so we never need to compare blocks 
$B_1, B_2$ such that $s\in B_1,B_2$.
\end{defn}

\medskip
It is not hard to check that both of the 
above block orders satisfy the lengthening condition.  
The lengthening condition and these
two block orders were introduced in [HK].

\begin{thm}\label{lampart}
Any block order satisfying the lengthening condition
yields a partitioning for $\Delta(\Pi_n^*)/S_{\lambda }$.  Hence,
$\Delta(\Pi_n)/S_{\lambda }$ is partitionable using the length-lex
order.
\end{thm}

\proof
Let us modify the chain-labeling for $\Delta (\Pi_n^*)/S_n$ as follows:
label covering relations with ordered 
4-tuples $(i,w_B,W,r)$ where $i$ is the number of bars to the left of the 
bar being inserted, 
$w_B$ is the content of the block immediately to the left of this bar, 
$W$ is the word obtained by concatenating all the block words 
to the left of the new bar and interspersing bar symbols between
the block words, and let $r$ be the rank at which the parent block $P$ 
into which the bar is inserted was itself created.  Precedence in the 4-tuple
proceeds from left to right.  The words $w_B$ are ordered by
a block order satisfying the lengthening condition.  
The words $W$ are ordered by considering the first block
where two words differ and then using our block order to compare these 
blocks.  In the
partitioning for $\Delta (\Pi_n^*)/S_{\lambda }$, 
we must use block content as well as size to determine block
equivalence, and we use our chosen block order that satisfies the lengthening 
condition to decide which offspring blocks are left children and which
are right children and also how to sort blocks, but otherwise the proof will
be identical to that in [He, p. 14-24], 
by virtue of the properties of the lengthening
discussion, to be discussed next.  

We claim that the lengthening condition ensures that we may 
replace any pair of consecutive bar insertions which proceed either (1) from 
right to left bar or (2) which insert bars from left to right in a single block
creating left children decreasing in size from left to right, by a 
lexicographically earlier saturated chain which overlaps ours in a 
codimension one face, yielding a topological descent in our saturated 
chain at the rank in between the two bar insertions.  In case 1, if the two
refinement steps refine distinct blocks from right to left, then it is clear
that these may also be refined left to right.  The other possibility for 
case 1 is 
that consecutive refinement steps split a single block $B$ into children 
$B_1,B_2,B_3$ by first refining $B$ into children $B_1B_2,B_3$, and that
$B_1B_2 < B_3$, $B_1<B_2$, so that bar insertions proceed right to left;
then the lengthening condition asserts that $B_1<B_1B_2$, so that the 
refinement first to children $B_1,B_2B_3$ gives a lexicographically smaller
chain.  In case 2, the first step refines a block $B$ into children
$B_1,B_2B_3$ and the next step refines $B_2B_3$
into children $B_2,B_3$.  We have that $B_1>B_2$ and $B_2<B_3$, which 
implies that $B_1 > \min (B_2,B_1B_3)$. 
We get a lexicographically smaller chain by first splitting 
$B$ instead into children $B_2,B_1B_3$.  These implications of the 
lengthening condition are adopted from [HK].
\EOP

\medskip
Now let us use the partitioning for $\Delta (\Pi_n)/S_1\times S_{n-1}$
derived from the
distinguished block order to obtain the following result.

\begin{thm}
If $S=\{ 1,\dots,i\} $, then $b'_S(n) = 1$ and otherwise $b'_S(n)>1$.
\end{thm}

\proof
Let $\{ s,t,\dots ,t\} $ denote the set of objects to be partitioned.
Notice that one may construct a saturated chain in $\Delta(\Pi_n^*)/S_{n-1}
\times S_1$ achieving any desired collection $S^*$ of topological descents by
inserting bars in the ordered set $s,t,\dots ,t$ 
as follows: for each maximal (possibly empty) string of 
ascents followed by a descent, we
place bars left to right filling the rightmost collection of available spots.
Finally, we insert bars left to right for the terminal string of ascents,
if there is one.  In this fashion, we achieve any $S$, implying 
$b'_S(n)\ge 1$ for all $S$.
For example, $S^* = \{ 1,4,5\} $ 
is achieved in $\Delta(\Pi_8^*)/S_1 \times S_7$
by $s|_6 t|_7 t|_5 t|_2 t|_3 t|_4 t|_ 1 t$ (letting subscripts denote ranks
of bar insertions) since this gives $w(S^*)=DA^2 D^2 A$.

Next we show $b'_{1,\dots ,i}(n)=1$.
Consider a word $w(S^*)=A^{n-i-2}D^i$, namely the case of 
$S=\{ 1,\dots i \} $.  A bar insertion isolating the $s$ followed by any
other bar insertion must comprise an ascent, while any bar insertion 
immediately before one isolating the unique $s$
must be a descent.  Thus, for $S=1,\dots ,i$, the $s$
must be isolated in either the first or the last refinement step.
It cannot be the first step since $b_{1,\dots ,j}(n) =0$, which means 
it would be 
impossible to refine the remaining nontrivial block of $n-1$ identical letters 
achieving a word $A^{n-j-3}D^i$.  Thus, the 
step splitting off the $s$ must come last.  Until the 
first descent, bars must be inserted   
left to right creating blocks of nondecreasing size (after the first block
which is automatically smallest by virtue of
containing the $s$).  The rightmost of these newly created 
blocks must have size 1, to avoid having a later ascent at any point
after the first descent.  Since the rightmost block has size one,
these increasing blocks all must have size one.  
Thus, only the block containing $s$ may be nontrivial after the initial 
series of consecutive ascents, so we must begin 
by inserting bars from left to right distance one 
apart filling up the rightmost available set of positions.
Now to avoid further ascents, we have no choice but to 
proceed right to left refining the
block containing $s$.
Since there is only one such saturated chain, we conclude
that $b'_{\{ 1,\dots ,i \} } (n)=1$.

For $S\ne \{ 1,\dots ,i \} $, we obtain a facet achieving $S$ as in the first
paragraph, but the fact that we have
a descent immediately before an ascent, gives enough
flexibility to guarantee an alternative facet also achieving $S$, constructed
as follows: the bar for the descent immediately preceding a string of ascents
may be placed one position farther to left
than the above greedy algorithm would choose.  If the
string of ascents is followed by a descent, then we put a bar at the 
rightmost position that is still vacant when we encounter the descent;
if the string of ascents concludes the entire string,
then we may place the bar insertion for the 
last ascent into this rightmost position.  In any case, $b'_S(n)>1$ for 
$S\ne \{ 1,\dots ,i\} $.
\EOP

\end{document}